# ON A CONJECTURE CONCERNING SOME
# AUTOMATIC CONTINUITY THEOREMS

## M. El Azhari

ABSTRACT. Let A and B be commutative locally convex algebras with unit. A is assumed to be a uniform topological algebra. Let $\Phi$ be an injective homomorphism from A to B. Under additional assumptions, we characterize the continuity of the homomorphism $\Phi^{-1}/\operatorname{Im}\Phi$ by the fact that the radical (or strong radical) of the closure of $\operatorname{Im}\Phi$ has only zero as a common point with $\operatorname{Im}\Phi$. This gives an answer to a conjecture concerning some automatic continuity theorems on uniform topological algebras.

1. INTRODUCTION. Let A and B be commutative locally convex algebras with unit. A is assumed to be a uniform topological algebra. Let $\Phi$ be an injective homomorphism from A to B. Under which conditions is $\Phi^{-1}/\operatorname{Im}\Phi$ continuous?
   Under additional assumptions such as:
(1) A is weakly regular and functionally continuous, B an lmc algebra, and $\overline{(\operatorname{Im}\Phi)}$ (the closure of $\operatorname{Im}\Phi$) is a semisimple Q-algebra; or
(2) A is weakly $\sigma^*$-compact-regular, and $\overline{(\operatorname{Im}\Phi)}$ is a strongly semisimple Q-algebra;
it is shown in [5] that $\Phi^{-1}/\operatorname{Im}\Phi$ is continuous, which improves earlier results by Bedaa, Bhatt and Oudadess ([2]).
   The following examples show that the hypothesis $\overline{(\operatorname{Im}\Phi)}$ is a Q-algebra in (1) and (2) cannot be omitted.

Example 1. Let $A = C[0,1]$ be the algebra of all complex continuous functions on the closed unit interval $[0,1]$. A is a uniform Banach algebra under the supnorm. Since $M(A)$ is homeomorphic to $[0,1]$, it follows that A is weakly regular. Consider $B = C[0,1]$. For any countable compact subset K of $[0,1]$, and $f \in B$, we put $p_K(f) = \sup\{|f(x)|, x \in K\}$. B is a complete uT-algebra under the system $(p_K)_K$. Consider $\Phi: A \to B$, $\Phi(f) = f$. Then $\overline{(\operatorname{Im}\Phi)} = B$ is semisimple but not a Q-algebra. Clearly $\Phi^{-1}/\operatorname{Im}\Phi$ is not continuous.

Example 2. Let $A = C_b(R)$ be the algebra of all complex continuous bounded functions on the real line. A is a uniform Banach algebra under the supnorm. A is weakly $\sigma^*$-compact-regular [2, Remark (4)]. Let $B = C(R)$ be the algebra of all complex continuous functions on R, with the compact-open topology. Consider $\Phi: A \to B$, $\Phi(f) = f$. Then $\overline{(\operatorname{Im}\Phi)} = C(R)$ is strongly semisimple but not a Q-algebra. Clearly $\Phi^{-1}/\operatorname{Im}\Phi$ is not continuous.

   In [2], the authors conjectured that the semisimplicity of $\overline{(\operatorname{Im}\Phi)}$ in (1) (and strong semisimplicity of $\overline{(\operatorname{Im}\Phi)}$ in (2)) can be omitted. According to the proofs in [2] and [5], the semisimplicity of $(Im\Phi)^-$ in (1) can be replaced by $Im\Phi \cap R((Im\Phi)^-) = \{0\}$, and the strong semisimplicity of $(Im\Phi)^-$ in (2) can be replaced by $Im\Phi \cap SR((Im\Phi)^-) = \{0\}$.

In this paper, we show that if A is weakly regular and functionally continuous, B an lmc algebra, and $(Im\Phi)^-$ is a Q-algebra, then the continuity of $\Phi^{-1}/Im\Phi$ is equivalent to $Im\Phi \cap R((Im\Phi)^- = \{0\}$. We also show that if A is weakly $\sigma^*$-compact-regular, B has continuous product, and $(Im\Phi)^-$ is a Q-algebra, then the continuity of $\Phi^{-1}/Im\Phi$ is equivalent to $Im\Phi \cap SR((Im\Phi)^- = \{0\}$.

2. PRELIMINARIES. All algebras considered are over the field $\mathbb{C}$, commutative, and having a unit element. A topological algebra is an algebra which is also a Hausdorff topological vector space such that the multiplication is separately continuous. A locally convex algebra (lc algebra) is a topological algebra whose topology is locally convex. A locally multiplicatively convex algebra (lmc algebra) is a topological algebra whose topology is determined by a family of submultiplicative seminorms. A uniform seminorm on an algebra A is a seminorm p such that $p(x^2) = p(x)^2$ for all $x \in A$. Such a seminorm is submultiplicative [4]. A uniform topological algebra (uT-algebra) is a topological algebra whose topology is determined by a family of uniform seminorms. A uniform normed algebra is a normed algebra $(A, \|\ \|)$ such that $\|x^2\| = \|x\|^2$ for all $x \in A$. Let A be an algebra and $x \in A$, we denote by $sp_A(x)$ the spectrum of x and $r_A(x)$ the spectral radius of x. For an algebra A, $M^*(A)$ denotes the set of all nonzero multiplicative linear functionals on A. For a topological algebra A, $M(A)$ denotes the set of all nonzero continuous multiplicative linear functionals on A. A topological algebra A is functionally continuous if $M^*(A) = M(A)$. A topological algebra is a Q-algebra [7] if the set of invertible elements is open. A topological algebra is weakly regular [2] if given a closed subset F of $M(A)$, $F \neq M(A)$, there exists a nonzero $x \in A$ such that $f(x) = 0$ for all $f \in F$. A topological algebra A is weakly $\sigma^*$-compact-regular [2] if given a compact subset K of $M^*(A)$, $K \neq M^*(A)$, there exists a nonzero $x \in A$ such that $f(x) = 0$ for all $f \in K$. We use $R(A)$ to denote the radical of an algebra A. If $R(A) = \{0\}$, we say that A is semisimple. Let A be a topological algebra with $M(A) \neq \emptyset$, the set $\{x \in A, f(x) = 0$ for all $f \in M(A)\}$ is called the strong radical of A and denoted by $SR(A)$. If $SR(A) = \{0\}$, we say that A is strongly semisimple. Let A be an lmc algebra, if A is complete or a Q-algebra, then $R(A) = SR(A)$.

3. RESULTS

Theorem 3.1. Let A be a weakly regular, functionally continuous, uT-algebra. Let B be an lmc algebra, and let $\Phi: A \to B$ be a one-to-one homomorphism such that $(Im\Phi)^-$ is a Q-algebra. Then the following are equivalent:
(1) $\Phi^{-1}/Im\Phi$ is continuous.
(2) $Im\Phi$ is functionally continuous.
(3) $\Phi^*: M((Im\Phi)^-) \to M(A)$, $\Phi^*(f) = f \circ \Phi$, is surjective.
(4) $Im\Phi \cap R((Im\Phi)^-) = \{0\}$.

Proof: (1) => (2): Let $F \in M^*(Im\Phi)$, $F = F \circ \Phi \circ (\Phi^{-1}/Im\Phi)$ is continuous since $F \circ \Phi$ and $\Phi^{-1}/Im\Phi$ are continuous.

(2) => (3): Let $f \in M(A)$ and $F = f \circ (\Phi^{-1}/\text{Im}\Phi)$, $F \in M^*(\text{Im}\Phi) = M(\text{Im}\Phi)$ and $f = F \circ \Phi$. $F$ can be extended to an $\bar{F} \in M((\text{Im}\Phi)^-)$. We have $f = \bar{F} \circ \Phi$. This shows that $\Phi^*$ is surjective.

(3) => (1): By [5, Theorem 2.1], the topology of A is defined by a family $\{p_s, s \in S\}$ of submultiplicative seminorms such that (i) for all $x \in A$ and $s \in S$ with $p_s(x) = 1$, there exists $f \in M(A)$ such that $|f(x)| = 1$. Let $s \in S$ and $y \in \text{Im}\Phi$ with $p_s(\Phi^{-1}(y)) \neq 0$. By (i), there exists $f \in M(A)$ such that $|f(\Phi^{-1}(y))| = p_s(\Phi^{-1}(y))$. Since $\Phi^*$ is surjective, there exists $F \in M((\text{Im}\Phi)^-)$ such that $f = F \circ \Phi$. We have $p_s(\Phi^{-1}(y)) = |f(\Phi^{-1}(y))| = |F(y)| \leq r_C(y)$, where $C = (\text{Im}\Phi)^-$. Since $C$ is a Q-algebra, $r_C$ is continuous at 0 [7, Proposition 13.5]. Then $\Phi^{-1}/\text{Im}\Phi$ is continuous.

(3) => (4): Let $y \in \text{Im}\Phi \cap R((\text{Im}\Phi)^-)$, there exists $x \in A$ such that $y = \Phi(x)$ and $F(\Phi(x)) = 0$ for all $F \in M((\text{Im}\Phi)^-)$. Then $f(x) = 0$ for all $f \in M(A)$ since $\Phi^*$ is surjective. Hence $x = 0$ and so $y = \Phi(x) = 0$ since $A$ is a uT-algebra.

(4) => (3): $\Phi^*$ is well defined and continuous. Since $(\text{Im}\Phi)^-$ is a Q-algebra, $M((\text{Im}\Phi)^-)$ is compact [6, p.187], thus $\Phi^*(M((\text{Im}\Phi)^-))$ is compact. Suppose that $\Phi^*$ is not surjective. By the weak regularity of A, there exists a nonzero $x \in A$ such that $f(\Phi(x)) = 0$ for all $f \in M((\text{Im}\Phi)^-)$. Since $(\text{Im}\Phi)^-$ is a Q-algebra, it follows that $\Phi(x) \in \text{Im}\Phi \cap R((\text{Im}\Phi)^-) = \{0\}$, and then $x = 0$, a contradiction.

Theorem 3.2. Let A be a weakly $\sigma^*$-compact-regular, uT-algebra. Let B be an lc algebra with continuous product, and $\Phi: A \to B$ be a one-to-one homomorphism such that $(\text{Im}\Phi)^-$ is a Q-algebra. The following are equivalent:
(1) $\Phi^{-1}/\text{Im}\Phi$ is continuous.
(2) $\text{Im}\Phi \cap SR((\text{Im}\Phi)^-) = \{0\}$.
(3) $\Phi^{**}: M((\text{Im}\Phi)^-) \to M^*(A)$, $\Phi^{**}(f) = f \circ \Phi$, is surjective.

Proof: (1) => (2): The topology of A is determined by a family $\{p_u, u \in U\}$ of uniform seminorms. For each $u \in U$, let $N_u = \{x \in A, p_u(x) = 0\}$ and $A_u$ be the Banach algebra obtained by completing $A/N_u$ in the norm $\|x_u\|_u = p_u(x)$, $x_u = x + N_u$. It is clear that $A_u$ is a uniform Banach algebra. For each $u \in U$, let $M_u(A) = \{f \in M(A), |f(x)| \leq p_u(x)$ for all $x \in A\}$. Let $u \in U$ and $x \in A$, $p_u(x) = \|x_u\|_u = r_u(x_u) = \sup\{|g(x_u)|, g \in M(A_u)\} = \sup\{|f(x)|, f \in M_u(A)\}$ by [7, Proposition 7.5] ($r_u$ is the spectral radius on $A_u$). Let $u \in U$ and $y \in \text{Im}\Phi$, $p_u(\Phi^{-1}(y)) = \sup\{|f \circ \Phi^{-1}(y))|, f \in M_u(A)\}$. Let $f \in M_u(A)$, $f \circ \Phi^{-1} \in M(\text{Im}\Phi) = M((\text{Im}\Phi)^-)$ since $\Phi^{-1}/\text{Im}\Phi$ is continuous and B has continuous product. Then $p_u(\Phi^{-1}(y)) \leq \sup\{|F(y)|, F \in M((\text{Im}\Phi)^-)\}$ for all $u \in U$ and $y \in \text{Im}\Phi$. Let $y \in \text{Im}\Phi \cap SR((\text{Im}\Phi)^-)$, we have $p_u(\Phi^{-1}(y)) = 0$ for all $u \in U$, then $\Phi^{-1}(y) = 0$ and so $y = 0$.

(2) => (3): $\Phi^{**}$ is continuous. Since $(\text{Im}\Phi)^-$ is a Q-algebra, $M((\text{Im}\Phi)^-)$ is compact [6, p.187], thus $\Phi^{**}(M((\text{Im}\Phi)^-))$ is compact. Suppose that $\Phi^{**}$ is not surjective. Since A is $\sigma^*$-compact-regular, there exists a nonzero $x \in A$ such that $f(\Phi(x)) = 0$ for all $f \in M((\text{Im}\Phi)^-)$. This gives $\Phi(x) \in \text{Im}\Phi \cap SR((\text{Im}\Phi)^-) = \{0\}$, and then $x = 0$, a contradiction.

(3) => (1): Similar to the proof of (3) => (1) in Theorem 3.1.

Here is an example such that A is weakly regular, uniform Banach algebra, B is a Banach algebra, $\Phi: A \to B$ is a one-to-one homomorphism, and $\text{Im}\Phi \cap R((\text{Im}\Phi)^-) = \{0\}$ but $R((\text{Im}\Phi)^-) \neq \{0\}$.

Example. Let $A = C[0,1]$ be the algebra of all complex continuous functions on the closed unit interval $[0,1]$. A is a uniform Banach algebra under the supnorm $\|.\|$, A is also weakly regular. By [3], there exists a norm $|.|$ on $C[0,1]$ such that $C[0,1]$ is an incomplete normed algebra. It is well known that $\|.\| \leq |.|$. Let B be the completion of $C[0,1]$ under the norm $|.|$. Consider $\Phi: A \to B$, $\Phi(f) = f$, we have $(Im\Phi)^- = B$. If B is semisimple, then $\Phi$ is continuous, and consequently the norms $\|.\|$ and $|.|$ are equivalent, a contradiction. Since $\|.\| \leq |.|$, $\Phi^{-1}/Im\Phi$ is continuous and so $Im\Phi \cap R((Im\Phi)^-) = \{0\}$ by Theorem 3.1.

Remark. The algebra A considered in the above example is also $\sigma^*$-compact-regular and $Im\Phi \cap SR((Im\Phi)^-) = \{0\}$ but $SR((Im\Phi)^-) \neq \{0\}$.

The following result is an application of Theorem 3.1.

Theorem 3.3. Let A be a functionally continuous normed algebra. Then the following assertions are equivalent:
(1) A is a uniform normed algebra.
(2) A has a largest closed, idempotent, absolutely convex, bounded subset.

Proof. (1) => (2): Let $\|.\|$ be a uniform norm defining the topology of A. Let $B = \{x \in A, \|x\| \leq 1\}$, B is a closed, idempotent, absolutely convex, bounded subset of A. Let C be an idempotent bounded subset of A. There exists $M > 0$ such that $\|x\| \leq M$ for all $x \in C$. Let $x \in C$, $\|x\| = \|x^{2^n}\|^{2^{-n}} \leq M^{2^{-n}}$ for all $n \geq 1$, then $\|x\| \leq 1$ i.e. $x \in B$.
(2) => (1): Let B be a largest closed, idempotent, absolutely convex, bounded subset of A. By [1, Proposition 2.15], we have $A = A(B) = \{tx, t \in \mathbb{C}$ and $x \in B\}$. Let $\|.\|_B$ be the Minkowski functional of B, $(A, \|.\|_B)$ is a normed algebra. By [1, Proposition 2.15], $\beta = \|.\|_B$ where $\beta$ is the radius of boundedness, then $(A, \|.\|_B)$ is a uniform algebra since $\beta(x^2) = \beta(x)^2$ for all $x \in A$. Let $A_1$ be the completion of A under the original norm. It is clear that $\Phi: (A, \|.\|_B) \to A_1$, $\Phi(x) = x$, is continuous, and consequently $(A, \|.\|_B)$ is functionally continuous. We now remark that we have proved the equivalence of (1), (2) and (3) in Theorem 3.1 without the condition that A is weakly regular. Using this remark, $\Phi^{-1}/Im\Phi$ is continuous, then $\Phi$ is a homeomorphism (into), so A is a uniform normed algebra.

Ecole Normale Supérieure
Avenue Oued Akreuch
Takaddoum, BP 5118, Rabat
Morocco

E-mail: mohammed.elazhari@yahoo.fr